\documentclass[letterpaper,12pt]{article}

\title{\bf The tropical Grassmannian}

\author{David Speyer \ and \ 
 Bernd Sturmfels \!\!\!\! \thanks{Partially supported by the
National Science Foundation (DMS-0200729)}
 \\
 \\
{\small Department of Mathematics,
  University of California, Berkeley} \\
{\small
{\tt $\{$speyer$,$bernd$\}$@math.berkeley.edu}
} }

\date{}

\usepackage{latexsym,array,delarray,amsthm,amssymb,psfrag}%eufrak
\usepackage[sumlimits]{amsmath}
\usepackage[dvips]{graphicx}
\theoremstyle{plain}
\newtheorem{thm}{Theorem}[section]
\newtheorem{lemma}[thm]{Lemma}
\newtheorem{prop}[thm]{Proposition}
\newtheorem{cor}[thm]{Corollary}

\newtheorem{rmk}[thm]{Remark}

\theoremstyle{definition}

\newtheorem{ex}[thm]{Example}

\theoremstyle{remark}

\newcommand{\zz}{\mathbb{Z}}
\newcommand{\nn}{\mathbb{N}}
\newcommand{\pp}{\mathbb{P}}

\newcommand{\rr}{\mathbb{R}}
\newcommand{\RR}{\mathbb{R}}

\newcommand{\cc}{\mathbb{C}}

\begin{document}
\maketitle

\begin{abstract}
\noindent
In tropical algebraic geometry, the solution sets of polynomial equations
are piecewise-linear. We introduce the tropical variety of a polynomial ideal,
and we identify it with a polyhedral subcomplex of the Gr\"obner fan.
The tropical Grassmannian arises in this manner from the
ideal of quadratic Pl\"ucker relations.
It  parametrizes  all tropical linear spaces.
Lines in tropical projective space are trees,
and their tropical Grassmannian ${\mathcal G}_{2,n}$ equals
the space of phylogenetic trees studied by
Billera, Holmes and Vogtmann. Higher
Grassmannians offer a natural generalization of the space of trees.
Their faces correspond to monomial-free initial ideals
of the Pl\"ucker ideal.
The tropical Grassmannian ${\mathcal G}_{3,6}$ is a 
simplicial complex glued from 1035 tetrahedra.
\end{abstract}

\section{Introduction}

The \emph{tropical semiring} 
$\,(\rr \, \cup \,\{\infty\}, {\rm min} ,+ )\,$
is the set of real numbers augmented by infinity
with the \emph{tropical addition},
which is taking the minimum of two numbers,
 and the \emph{tropical multiplication}
which is the ordinary addition \cite{Pin}.
These operations satisfy the familiar
axioms of arithmetic, e.g.~distributivity,
with $\infty$ and $0$ being the two neutral elements.
\emph{Tropical monomials}
$\, x_1^{a_1}  \cdots x_n^{a_n}\,$
represent ordinary linear forms $\sum_i a_i x_i$,
and \emph{tropical polynomials}
\begin{equation}
\label{tropicalpol}
 F(x_1,x_2,\ldots,x_n) \quad = \quad
\sum_{ a \in {\mathcal A}}  C_a\, x_1^{a_1} x_2^{a_2} \cdots x_n^{a_n},
\,\quad \hbox{with}\,\,{\mathcal A} \subset \nn^n, 
C_a \in \rr,
\end{equation}
represent piecewise-linear convex functions
$\ F : \rr^n \rightarrow \rr $.
To compute $F(x)$, we take the minimum of the 
affine-linear forms  $\,C_a + \sum_{i=1}^n a_i x_i \,$
for  $\,a \in {\mathcal A} $.
We define the \emph{tropical hypersurface} ${\mathcal T}(F)$
as the set of all points $x $ in $\rr^n$ for which 
this minimum is attained at least twice, as $a$ runs over $
{\mathcal A}$. Equivalently, ${\mathcal T}(F)$
is the set of all points $x \in \rr^n$ at which $F$ is not differentiable.
Thus a tropical hypersurface is an $(n-1)$-dimensional
polyhedral complex  in $\rr^n$.

The rationale behind this definition will become clear 
in Section 2, which gives a self-contained
development of the basic theory of  tropical varieties.
For further background and pictures see \cite[\S 9]{SSPE}.
Every  tropical variety is a finite intersection of
tropical hypersurfaces (Corollary \ref{FiTG}).
But not every intersection of tropical hypersurfaces
is a tropical variety (Proposition \ref{caterpillar}).
Tropical varieties are also known as  \emph{logarithmic limit sets}
\cite{Berg}, \emph{Bieri-Groves sets} \cite{BG}, or
 \emph{non-archimedean amoebas} \cite{EKLW}.
Tropical curves are the key ingredient in
Mikhalkin's formula \cite{Mik} for
planar Gromov-Witten invariants.

The object of study in this paper is
the tropical Grassmannian ${\mathcal G}_{d,n}$
which is a polyhedral fan in $\rr^{\binom{n}{d}}$
defined by the ideal of quadratic Pl\"ucker relations.
All of our main results 
regarding ${\mathcal G}_{d,n}$ are stated in
Section 3. The proofs appear in the subsequent sections.
In Section 4 we prove Theorem \ref{phylogenetic} which identifies
${\mathcal G}_{2,n}$ with the
space of phylogenetic trees in \cite{BHV}.
A detailed study of the fan
${\mathcal G}_{3,6} \subset \rr^{20}$ is presented in Section 5. 
In Section 6 we introduce tropical linear spaces 
and we prove that they are parametrized by
the tropical Grassmannian (Theorem \ref{G36}).
In Section $7$ we show that the tropical Grassmannian
${\mathcal G}_{3,7}$ depends on the characteristic of
the ground field.

\section{The tropical variety of a polynomial ideal}

Let $K$ be an algebraically closed field with
a valuation into the reals, denoted
$\,{\rm deg} : K^* \rightarrow \rr $. 
We assume that 1 lies in the image of $\deg$ and we fix $t \in K^*$ with $\deg(t)=1$. The corresponding
local ring and its maximal ideal are
$$ 
R_K \,\,\, = \,\,\, \{ \,c \in K \,\,: \,\, {\rm deg}(c) \geq 0 \,\}
\quad \hbox{and} \quad
M_K \,\,\, = \,\,\, \{ \,c \in K \,\,: \,\, {\rm deg}(c) > 0 \,\}.  $$
The residue field $\, k \, = \,R_K/M_K \,$ is algebraically closed.
Given any ideal
$$ I \quad \subset \quad
\,K[x] \,\,\, = \,\,\, K[x_1,x_2,\ldots,x_n] , $$
we consider its affine variety in the
\emph{$n$-dimensional algebraic torus} over $K$,
$$ {V}(I) \quad = \quad
\bigl\{\, u \in (K^*)^n \,\, : \,\,
f(u)  \, = \, 0 \,\,\,\hbox{for all} \,\,
f \in I \,\bigr\}. $$
Here $K^* = K \backslash \{0\}$.
In all our examples, $K$ is the algebraic 
closure of the rational function field $\cc(t)$
and  ``{\rm deg}'' is the standard valuation
at the origin. Then $k = \cc$, and if $c \in \cc[t]$ then
$\,{\rm deg}(c) \,$ is the order of vanishing
of $c$ at $0$.
 
Every  polynomial in $K[x]$ maps to a 
tropical polynomial as follows. If
\begin{equation}
\label{regularpol}
f(x_1,\ldots,x_n) \,\,\, = \,\,\,
\sum_{ a \in {\mathcal A}}  c_a\, x_1^{a_1} \cdots x_n^{a_n}\,
\qquad \hbox{with $\, c_a \in K^* \,$ for
$a \in {\mathcal A}$}.
\end{equation}
and  $\,C_a = {\rm deg}(c_a)$, then 
$ {\rm trop}(f)$ denotes the tropical polynomial $F$
in (\ref{tropicalpol}).

The following definitions are a variation on
Gr\"obner basis theory \cite{GB+CP}.
Fix $w \in \rr^n$.  The \emph{$w$-weight} of a term
$\,c_a \cdot x_1^{a_1}\cdots x_n^{a_n} \,$ in
(\ref{regularpol}) is
$\,{\rm deg}(c_a) + a_1 w_1 + \cdots + a_n w_n$.
The  \emph{initial form} $\,{\rm in}_w(f)\,$ 
of a polynomial $f$ is defined as follows. 
Set $\tilde{f}(x_1, \ldots, x_n)=f(t^{w_1} x_1, \ldots, t^{w_n} x_n)$.
 Let $\nu$ be the smallest weight of any term of $f$, 
so that $t^{-\nu} \tilde{f} $ is a non-zero element in
$\,R_K[x]$.  Define ${\rm in}_w(f)$ 
as the image of $t^{-\nu} \tilde{f}$ in $k[x]$. 
We set ${\rm in}_w(0)=0$. For
$K=\overline{\cc(t)}$ and $k=\cc$ this means that the initial form
 ${\rm in}_w(f) $ is a polynomial in  $\cc[x]$.

Given any ideal $I \subset K[x]$, then its \emph{initial ideal} is
defined to be
$$  {\rm in}_w(I) \quad = \quad
\bigl\langle \, {\rm in}_w(f) \,\, : \,\,
f \in I \,\,  \bigr\rangle
\quad \subset \quad k[x] . $$

\begin{thm}
\label{charac}
For an ideal $I \subset K[x]$ the following 
subsets of $\,\rr^n$ coincide:
\begin{itemize}
\item[(a)] The closure of the set $\, \bigl\{
( {\rm deg}(u_1),\ldots,{\rm deg}(u_n))
\,: \, (u_1,\ldots,u_n) \in {V}(I) \bigr\}$;
\item[(b)] The intersection of the tropical hypersurfaces
${\mathcal T}({\rm trop}(f))$ where 
$f \in I $;
\item[(c)] The set of all vectors $w \in \rr^n$ such that
${\rm in}_w(I)$ contains no monomial.
\end{itemize}
\end{thm}

The set defined by the three conditions
in Theorem \ref{charac} is denoted ${\mathcal T}(I)$
and is called the \emph{tropical variety} of the ideal $I$.
Variants of this theorem already appeared in 
\cite[Theorem 9.17]{SSPE} and in \cite[Theorem 6.1]{EKLW},
without and with proof respectively.
Here we present a short proof which is self-contained.

\begin{proof}
First consider any point
$\,w = ( {\rm deg}(u_1),\ldots,{\rm deg}(u_n))\,$
in the set (a). For any $f \in I$ we have
$\,f(u_1,\ldots,u_n) = 0 \,$ and this implies
that the minimum in the definition of
$F = {\rm trop}(f)$ is attained at least
twice at $w$. This condition is equivalent to
${\rm in}_w(f)$ not being a monomial.
This shows that (a) is contained in (b), and
(b) is contained in (c). It remains to prove
that (c) is contained in (a).  Consider any vector
$w$ in (c) such that
$w =  ( {\rm deg}(v_1),\ldots,{\rm deg}(v_n))\,$
for some $v \in (K^*)^n$. Since
the image of the valuation is dense in $\rr$
and the set defined in (a) is closed,
it suffices to prove that
$\,w = ( {\rm deg}(u_1),\ldots,{\rm deg}(u_n))$
for some $u \in {V}(I)$. By making the
change of coordinates $ \, x_i \, = \,x_i \cdot v_i^{-1}$,
we may assume that $\, w = (0,0,\ldots,0)$.

Since ${\rm in}_w(I)$ contains no monomial and
since $k$ is algebraically closed, by
the Nullstellensatz there exists a point
$\,\bar{u} \in {V}({\rm in}_w(I)) \subset (k^*)^n$.
Let $\bar{m}$ denote the maximal ideal in $k[x]$
corresponding to $\bar{u}$.
Let $S$ be the set of polynomials $f$ in
 $\,R_K[x]\,$  whose reduction modulo $M_K$
is not in $\bar{m}$. Then
$S$ is a  multiplicative set in $R_K[x]$ 
disjoint from $I$. Consider the induced map
$$\phi \,\,:\,\, R_K \,\,\, \longrightarrow 
\,\,\, S^{-1} R_K[x]/S^{-1} (I \cap R_K[x])$$
Let $P$ be a minimal prime of the ring on the right
hand side.  We claim that $\phi^{-1}(P) = \{0\}$.
Suppose not, and pick $c \in R_K \backslash \{0\}$
with $\phi(c) \in P$. Then $\phi(c)$ is 
a zero-divisor in  $\,R_K[x]/S^{-1} (I \cap R_K[x])$,
so we can find $f \in S$ such that  $cf \in I$.
Since $c^{-1}$ exists in $K$, this implies
$f \in I$ which is a contradiction.

Now, $\phi^{-1}(P) = \{0\}$ implies that
$\,P \otimes_{R_K} K \,$ is a proper ideal in $K[x]/I$.
There exists a maximal ideal
of $K[x]/I$ containing $\,P \otimes_{R_K} K $,
and, since $K$ is algebraically closed, this
maximal ideal has the form
$\,  \langle x_1 - u_1,\ldots,x_n-u_n \rangle\,$
for some $u \in {V}(I) \subset (K^*)^n$.
We claim that $u_i \in R_K$ and
$ u_i \cong \bar{u}_i \mod M_K$.
This will imply ${\rm deg}(u_1)= \cdots =  {\rm deg}(u_1) = 0$
and hence complete the proof.

Consider any $x_i-u_i \in I$.
By clearing denominators, we get 
$a_i x-b_i \in I \cap R_K[x]$ with $b_i/a_i=u_i$,
and not both $a_i$ and $b_i$ lie in $M_k$. 
If $a_i \in M_K$, then $a_i x-b_i \cong -b_i \mod M_K$.
Hence ${\rm in}_w(I)$ contains $b_i \in K^*$
and hence equals the unit ideal, which is a contradiction.
If $a_i \not\in M_K$ and $-b_i/a_i \not\cong \bar{u}_i \mod M_K$ 
then the reduction of $a_i x -b_i$ modulo $M_K$ 
does not lie in $\bar{m}$. This means that
$a_i x -b_i \in S$ and is a unit of $S^{-1}R_K[x]$, 
so $P$ is the unit ideal. But then $P$ is not prime,
also a contradiction. This completes the proof.
\end{proof}

The key point in the previous proof can be summarized as follows:

\begin{cor}
Every zero over $k$ of the initial ideal ${\rm in}_w(I)$ 
lifts to a zero over $K$ of $I$.
\end{cor}

By zero of an ideal $I$ in $K[x]$ we mean
a point on its variety in $(K^*)^n$. 
The notion of (reduced) Gr\"obner bases
is well-defined for ideals  $I \subset K[x]$
and (generic) weight vectors $w$, and, 
by adapting the methods of \cite[\S 3]{GB+CP}
to our situation, we can 
compute a \emph{universal Gr\"obner basis}   ${\rm UGB}(I)$.
This is a  finite subset of $I$
which contains a Gr\"obner basis for $I$
with respect to any weight vector  $w \in \rr^n$.
{}From part (c) of  Theorem \ref{charac} we derive:

\begin{cor}  \label{FiTG}
{\rm [Finiteness in Tropical Geometry]}
The tropical variety ${\mathcal T}(I)$ is the
intersection of  the tropical hypersurfaces
${\mathcal T}({\rm trop}(f))$ where 
$f \in {\rm UGB}(I)$.
\end{cor}

The following result is due to Bieri and Groves \cite{BG}.
An alternative proof using Gr\"obner basis methods appears in
\cite[Theorem 9.6]{SSPE}.

\begin{thm} {\rm [Bieri-Groves Theorem]}
\label{bierigroves}
If $I$ is a prime ideal and $K[x]/I$ has
Krull dimension $r$, then ${\mathcal T}(I)$ is 
a pure polyhedral complex of dimension $r$.
\end{thm}

We shall be primarily interested in the
case when $k = \cc$ and $K = \overline{\cc(t)}$.
Under this hypothesis, 
the ideal $I$ is said to have \emph{constant coefficients}
if the coefficients $c_a$ of the generators $f$ of $I$ lie
in the ground field $\cc$. This implies
$C_a = {\rm deg}(c_a) = 0$ in (\ref{tropicalpol}), where
$F = {\rm trop}(f)$. Our problem is now
 to solve a system of tropical
equations all of whose coefficients are 
identically zero:
\begin{equation}
\label{tropicalpolzero}
 F(x_1,x_2,\ldots,x_n) \quad = \quad
\sum_{ a \in {\mathcal A}} \,0 \cdot x_1^{a_1} x_2^{a_2} \cdots x_n^{a_n}.
\end{equation}
Here the tropical variety is a subfan of the
Gr\"obner fan of an ideal in $\cc[x]$. 

\begin{cor}
\label{fancoro}
If $I$ has constant coefficients then
${\mathcal T}(I)$ is a fan in $\rr^n$.
\end{cor}

\section{Results on the tropical Grassmannian}

We fix a polynomial ring in  $\binom{n}{d}$
variables with integer coefficients:
$$ \zz[p] \quad = \quad \zz\bigl[p_{i_1 i_2 \cdots i_d}
\,:\,1 \leq i_1 < i_2 < \cdots < i_d \leq n \bigr]. $$
The \emph{Pl\"ucker ideal} $I_{d,n}$ is the 
homogeneous prime ideal in $\zz[p]$ consisting 
of the algebraic relations among the
$d \times d$-subdeterminants of any
$d \times n$-matrix with entries in any
 commutative ring.
The ideal $I_{d,n}$ is generated by quadrics, and 
it has a well-known quadratic Gr\"obner basis
(see e.g.~\cite[Theorem 3.1.7]{AIT}).
The projective variety of $I_{d,n}$ is
the \emph{Grassmannian} $G_{d,n}$ which
parametrizes all $d$-dimensional linear subspaces
of an $n$-dimensional vector space.

The \emph{tropical Grassmannian} ${\mathcal G}_{d,n}$
is the tropical variety ${\mathcal T}(I_{d,n})$
of the Pl\"ucker ideal $I_{d,n}$,
over a field $K$ as in Section 2.
 Theorem \ref{charac} (c) implies
$$ {\mathcal G}_{d,n} \quad = \quad
\bigl\{ \, w \in \rr^{\binom{n}{d}} \,\,: \,\,
{\rm in}_w(I_{d,n}) \,\,\, 
\hbox{contains no monomial} \,\bigr\}. $$
The ring $(\zz[p]/I_{d,n}) \otimes K$ 
is known to have Krull dimension $(n-d)d+1$.
Therefore Theorem \ref{bierigroves} and Corollary \ref{fancoro} imply
the following statement.

\begin{cor} The tropical Grassmannian   $\,{\mathcal G}_{d,n}\,$ is
a polyhedral fan in $\, \rr^{\binom{n}{d}} $. Each of its maximal
cones has the same dimension, namely, $(n-d)d+1$.
\end{cor}

We show in  Section 7 that
the fan ${\mathcal G}_{d,n}$ depends on the
characteristic of $K$ if $d \geq 3$ and $n \geq 7 $.
All results in Sections 2--6 are valid over
any field $K$.

It is convenient to reduce the dimension of
the tropical Grassmannian. This can be done in
three possible ways.  Let $\phi$ denote the linear map
from $\,\rr^n \,$ into $\,\rr^{\binom{n}{d}}\,$
which sends an $n$-vector $(a_1,a_2,\ldots,a_n)$ 
to the $\binom{n}{d}$-vector whose $(i_1, \ldots,i_d)$-coordinate
is $\,a_{i_1}  + \cdots + a_{i_d}$. The map $\phi$ is injective, 
and its image is the common intersection of all cones in the tropical
Grassmannian ${\mathcal G}_{d,n}$. Note that vector
$(1,\ldots,1)$ of length $\binom{n}{d}$ lies in ${\rm image}(\phi)$.
We conclude:

\begin{itemize}
\item
The image of  ${\mathcal G}_{d,n}$ in $\rr^{\binom{n}{d}}/ \rr (1,\ldots,1)$
is a fan ${\mathcal G}_{d,n}'$ of dimension $d(n-d)$.
\item
The image of ${\mathcal G}_{d,n}$
or ${\mathcal G}_{d,n}'$ in $\rr^{\binom{n}{d}}/{\rm image}(\phi)$
is a fan ${\mathcal G}_{d,n}''$ of dimension $(d-1)(n-d-1)$.
No cone in this fan contains a non-zero linear space.
\item Intersecting  ${\mathcal G}_{d,n}''$
with the unit sphere yields a polyhedral
complex ${\mathcal G}_{d,n}'''$. Each maximal face
of ${\mathcal G}_{d,n}'''$ is a polytope
of dimension $\,nd - n - d^2$.
\end{itemize}

We shall distinguish the four objects ${\mathcal G}_{d,n}$,
${\mathcal G}_{d,n}'$, ${\mathcal G}_{d,n}''$ and  
${\mathcal G}_{d,n}'''$  when stating our theorems below.
In subsequent sections less precision is needed, and we
sometimes identify ${\mathcal G}_{d,n}$,
${\mathcal G}_{d,n}'$, ${\mathcal G}_{d,n}''$ and  
${\mathcal G}_{d,n}'''$ if there is no danger of confusion.

\begin{ex}
\label{threepoints} $(d=2,n=4)$
The smallest non-zero Pl\"ucker ideal is
the principal ideal
$\, I_{2,4} \, = \, \langle p_{12} p_{34} - p_{13} p_{24} +
p_{14} p_{23} \rangle $. 
Here ${\mathcal G}_{2,4}\,$ is a fan
with three five-dimensional cones
$\,\rr^4 \times \rr_{\geq 0} \,$ glued
along $\,\rr^4 = {\rm image}(\phi)$.
The  fan  ${\mathcal G}_{2,4}''$
consists of three half rays emanating from the origin
(the picture of a tropical line).
The zero-dimensional simplicial complex
 ${\mathcal G}_{2,4}'''$
consists of \emph{three points}.
\end{ex}

\begin{ex}
\label{petersen} $(d =  2,n =  5)$
The tropical Grassmannian
${\mathcal G}_{2,5}'''$
is the \emph{Petersen graph} with
$10$ vertices and $15 $ edges.
This was shown in \cite[Example 9.10]{SSPE}.
\end{ex}

The following theorem generalizes both of
these examples. It concerns the case $d=2$,
that is, the tropical Grassmannian of
lines in $(n-1)$-space.

\begin{thm} \label{phylogenetic} $\!\!\!$ The tropical Grassmannian 
$ {\mathcal G}_{2,n}''' \!$ is a simplicial complex known as 
space of phylogenetic trees. It has $2^{n-1} - n-1 $ vertices, 
 $ 1 \cdot 3  \cdots (2n  \! - \!  5)$ facets, and its
homotopy type is a bouquet of
$\,(n \!- \!2) \, !$ spheres of dimension $n\! - \! 4$.
\end{thm}

A detailed description of $ {\mathcal G}_{2,n}$
and the proof of this theorem  will be given in Section 4.
Metric properties of the \emph{space of phylogenetic trees}
were studied by Billera, Holmes and Vogtmann in \cite{BHV} 
(our $n$ corresponds to Billera, Holmes and Vogtmann's $n+1$.)
The abstract simplicial complex and its homotopy type had 
been found earlier by Vogtmann \cite{V} and by Robinson and 
Whitehouse \cite{RW}. The description has the following
corollary. Recall that a simplicial complex is a  \emph{flag
complex} if the minimal non-faces are pairs of vertices.
This property is crucial for the existence of
unique geodesics in \cite{BHV}.

\begin{cor}
\label{isflagc}
The simplicial complex
$\, {\mathcal G}_{2,n}''' \,$
is  a flag complex.
\end{cor}

We do not have a complete description of the
tropical Grassmannian in the general case
$d \geq 3$ and $n-d \geq 3$.
We did succeed, however, in computing all monomial-free
initial ideals $\, {\rm in}_w(I_{d,n})\,$ for $d = 3$ and $n=6$:

\begin{thm}
\label{G36}
The tropical Grassmannian  ${\mathcal G}_{3,6}'''$
is a $3$-dimensional simplicial complex 
with $65$ vertices,
$550$ edges, $1395$ triangles and  $1035$ tetrahedra.
\end{thm}

The proof and complete description of 
${\mathcal G}_{3,6}$ will be presented in Section 5.
We shall see that ${\mathcal G}_{3,6}$ differs in
various ways from the tree space ${\mathcal G}_{2,n}$.
Here is one instance of this, which follows from Theorem~\ref{G36Detail}. Another one
is Corollary \ref{inibinoisprime}  versus Proposition 
\ref{radicalnotprime}.

\begin{cor}
\label{G36notflag}
The tropical Grassmannian  ${\mathcal G}_{3,6}'''$
is not a flag complex.
\end{cor}

If $X$ is  a $d$-dimensional linear subspace
of the vector space $K^n$, then (the
topological closure of) its image
${\rm deg}(X)$ under the degree map is 
a polyhedral complex in $\rr^n$. 
Such a polyhedral complex arising from a $d$-plane
in $K^n$ is called a  \emph{tropical $d$-plane in $n$-space}.
Since $X$ is invariant under scaling, every cone in 
${\rm deg}(X)$ contains the line spanned by
$(1,1,\ldots,1)$, so we can identify  ${\rm deg}(X)$ with 
its image in  $\, \rr^n/\rr (1,1,\ldots,1)
\simeq \rr^{n-1}$. Thus
${\rm deg}(X)$ becomes a  $(d-1)$-dimensional polyhedral complex
in $\,\rr^{n-1}$. For $d=2$, we get a tree.

The classical Grassmannian $G_{d,n}$ is the
projective variety in  $ \pp^{\binom{n}{d}-1}$
defined by the Pl\"ucker ideal $I_{d,n}$. There is a
canonical bijection between  $G_{d,n}$ and the
set of $d$-planes through the origin in $K^n$.
  The analogous bijection
for the tropical Grassmannian $ {\mathcal G}_{d,n}'$
is the content of the next theorem.

\begin{thm}
\label{grassmain}
 The bijection between the classical Grassmannian 
 $G_{d,n}$ and the set of $d$-planes in $K^n$ induces
a unique bijection  $\, w \mapsto L_w \,$
between the tropical Grassmannian
$\, {\mathcal G}_{d,n}'\,$
and the set of tropical $d$-planes in $n$-space.
\end{thm}

Theorems  \ref{phylogenetic}, 
\ref{G36} and \ref{grassmain}
are proved in Sections
4, 5 and 6 respectively.

\section{The space of phylogenetic trees}

In this section we prove  Theorem \ref{phylogenetic}
which asserts that  the  \emph{tropical Grassmannian of lines}
$\, {\mathcal G}_{2,n}\,$ coincides 
with the \emph{space of phylogenetic trees}
 \cite{BHV}. We begin by reviewing the
simplicial complex $\,{\bf T}_n\,$
underlying this space.

The vertex set ${\rm Vert}({\bf T}_n)$ consists of all unordered pairs
$\,\{A,B\}$, where $A$ and $B$ are disjoint subsets of 
$\,[n] := \{1,2,\ldots,n\}\,$
having cardinality at least two, and 
$\, A \,\cup \,B \, = \, [n]$. The 
cardinality of  ${\rm Vert}({\bf T}_n)$  is $\, 2^{n-1} - n - 1 $.
Two vertices $\,\{A,B\}\,$ and $\,\{A',B'\}$
are connected by an edge in ${\bf T}_n$ if and only if
\begin{equation}
\label{disjunction4}
  A \subset A'  \quad {\rm or} \quad  
  A \subset B' \quad {\rm or} \quad  
 B \subset A'  \quad {\rm or} \quad   B \subset B'.
\end{equation}
We now define ${\bf T}_n$ as the
flag complex with this graph. Equivalently,
a subset $\sigma \subseteq {\rm Vert}({\bf T}_n)$
is a face of ${\bf T}_n$ if any pair
 $\bigl\{ \{A,B\},\{A',B'\} \bigr\} \subseteq \sigma$
satisfies (\ref{disjunction4}).

The simplicial complex ${\bf T}_n$ was first introduced
by Buneman (see \cite[\S 5.1.4]{BaGu}) and was studied
more recently by Robinson-Whitehouse \cite{RW} and Vogtmann \cite{V}.
These authors obtained the following results.
Each face $\sigma$ of ${\bf T}_n$ corresponds
to a semi-labeled tree with leaves $1,2,\ldots,n$.
Here each internal node is unlabeled and has
at least three neighbors. Each internal edge
of such a tree defines a partition  $\{A,B\}$ 
of the set of leaves $\{1,2,\ldots,n\}$, and
we encode  the tree by the set of partitions
representing its internal edges.
The facets (= maximal faces) of ${\bf T}_n$ correspond to
\emph{trivalent trees}, that is, 
semi-labeled trees  whose internal nodes
all have three neighbors. All facets
of ${\bf T}_n$ have the same cardinality $n-3$,
the number of internal edges of any trivalent tree.
Hence ${\bf T}_n$ is pure of dimension $n-4$. The number of
facets (i.e.~trivalent semi-labeled trees
on $\{1,2,\ldots,n\}$) is
the \emph{Schr\"oder number}
\begin{equation}
\label{schroder}
 (2n-5) !! \quad = \quad
(2n-5) \times  (2n-7) \times \cdots \times 5 \times 3 \times 1. 
\end{equation}
It is proved in \cite{RW} and
\cite{V} that ${\bf T}_n$ has the homotopy type of a bouquet of
$\,(n-2) \, !\,$ spheres of dimension $\,n-4$.
The two smallest cases $n=4$ and $n=5$ are
discussed in Examples \ref{threepoints} and \ref{petersen}.
Here is a description of the next case:

\begin{ex} {\rm $(n=6)$ }
The two-dimensional simplicial complex ${\bf T}_6$ has $25$
vertices, $105$ edges and $105$ triangles, each coming
in two symmetry classes:
\begin{eqnarray*}
& \hbox{$15$ vertices like  }
\{12,3456\} \,,
\quad
 \hbox{$10$ vertices like  }
\{123,456\} , \\
& \hbox{$60$ edges like  } 
\{\{12, 3456\}, \{123, 456\}\} , \\
&  \hbox{$45$ edges like  } 
\{\{12, 3456\}, \{ 1234, 56 \}\}, \\
&  \hbox{$90$ triangles like  } 
\{\{ 12, 3456\},
 \{ 123, 456 \},
 \{ 1234, 56 \}\}, \\
&  \hbox{$15$ triangles like  } 
\{\{12, 3456\}\},
 \{ 34, 1256\}\},
 \{ 56 , 1234 \}\} .
\end{eqnarray*}
Each edge lies in three triangles,
corresponding to restructuring subtrees. \qed
\end{ex}

We next describe an embedding of ${\bf T}_n$ as a 
simplicial fan into the $\frac{1}{2}n(n-3)$-dimensional
vector space $\rr^{\binom{n}{2}}/{\rm image}(\phi)$.
For each trivalent tree $\sigma$
we  first define a cone $B_\sigma$ in 
$\rr^{\binom{n}{2}}$ as follows.
By a \emph{realization} of a semi-labeled tree $\sigma$
we mean a one-dimensional cell complex
in some Euclidean space whose underlying
graph is a tree isomorphic to $\sigma$.
Such a realization of $\sigma$ is a metric
space on $\{1,2,\ldots,n\}$.
The \emph{distance} between $i$ and $j$ is
the length of the unique path between leaf $i$
and leaf $j$ in that realization. Then we set
\begin{eqnarray*}
  B_\sigma \quad = \,\,\ &
\bigl\{  \, (w_{12}, w_{13}, \ldots, w_{n-1,n}) 
\in \rr^{\binom{n}{2}} \,\,:\,\,
-  w_{ij}  \,\,\hbox{
is the distance from} \\ & \hbox{leaf $i$ to  leaf $j$ in some
realization of $\sigma$} \bigr\}
\,\,\,\, + \,\,\, {\rm  image}(\phi).
\end{eqnarray*}
Let $C_\sigma$  denote the
image of  $  B_\sigma $ in the  quotient space
$\rr^{\binom{n}{2}}/{\rm image}(\phi)$.
Passing to this quotient has the geometric meaning
that two trees are identified if their only difference
is in the lengths of the $n$ edges adjacent to the leaves.

\begin{thm}
\label{bhsfan} $\!$
The closure $\overline{C}_\sigma$ is
a simplicial cone of dimension $|\sigma|$
with relative interior $C_\sigma$.
The collection of all cones $\overline{C}_\sigma$,
as $\sigma$ runs over ${\bf T}_n$,
is a simplicial fan. It is isometric to the
Billera-Holmes-Vogtmann space of trees.
\end{thm}

\begin{proof} 
Realizations of semi-labeled trees are characterized
by the \emph{four point condition}
(e.g.~\cite[Theorem 2.1]{BaGu}, \cite{Bu}).
This condition states
that for any quadruple of leaves $i,j,k,l$ there exists a unique
relabeling such that
\begin{equation}
\label{fourpoint}
w_{ij} + w_{kl} \,\, = \,\,
w_{ik} + w_{jl} \,\, \leq \,\,
w_{il} + w_{jk} . 
\end{equation}
Given any tree $\sigma$, this gives
a system of  $\binom{n}{4}$ linear equations and 
$\binom{n}{4}$ linear inequalities.
The solution set of this linear system
is precisely the closure $ \overline{B}_\sigma $
of the cone $B_\sigma$ in $\rr^{\binom{n}{2}}$.
This follows from the
\emph{Additive Linkage Algorithm}  \cite{Bu}  which 
reconstructs the combinatorial tree $\sigma$ from any
point $w$ in $B_\sigma$. 

All of our cones share a common linear subspace, namely,
\begin{equation}
\label{lineality}
 \overline{B}_\sigma \,\, \cap \,\, - \overline{B}_\sigma
\quad = \quad {\rm image}(\phi). 
\end{equation}
This is seen by replacing the inequalities in (\ref{fourpoint}) by equalities.
The cone $\overline{B}_\sigma$ is the direct sum (\ref{conegenerated})
of this linear space with a $|\sigma|$-dimensional simplicial cone.
 Let $\,\{ e_{ij} \, : \, 1 \leq i < j \leq n \}\,$ 
denote the standard basis of $\rr^{\binom{n}{2}}$.
Adopting the convention $e_{ji} =  e_{ij}$,
for any partition $\{A,B\}$ of $\{1,2,\ldots,n\}$ we define
$$ E_{A,B} \quad = \quad \sum_{i \in A} \sum_{j \in B} e_{ij} . $$
These vectors give the generators of our cone as follows:
\begin{equation}
\label{conegenerated}
\overline{B}_\sigma
\quad = \quad
 {\rm image}(\phi)\,\,\,\, + \,\,\,\,  \rr_{\geq 0}  \,
\bigl\{\, E_{A,B} \,: \,\{A,B\} \in \sigma \,\bigr\}. 
\end{equation}
{}From the two presentations (\ref{fourpoint})
and (\ref{conegenerated}) it follows that
\begin{equation}
\label{itsafan}
\overline{B}_\sigma \,\,\, \cap \,\, \overline{B}_\tau 
\quad = \quad
\overline{B}_{\sigma \,\cap\, \tau}
\qquad \hbox{for all} \,\, 
 \sigma,\tau \in {\bf T}_n. 
\end{equation}
Therefore the cones $B_\sigma$ form a fan
in $\rr^{\binom{n}{2}}$,
and this fan has face poset ${\bf T}_n$.
It follows from (\ref{conegenerated}) that the  
quotient $\,\overline{C}_\sigma 
\, = \,\overline{B}_\sigma /{\rm image}(\phi)\,$ is a pointed cone.

We get the desired conclusion for the cones
$\overline{C}_\sigma$  by taking quotients modulo the common
linear subspace (\ref{lineality}).
The resulting fan in 
$\rr^{\binom{n}{2}}/{\rm image}(\phi)$
is simplicial of pure dimension $n-3$
and has face poset ${\bf T}_n$.  It is
isometric to the Billera-Holmes-Vogtmann
space in \cite{BHV} because their
metric is flat on each cone
$\,\overline{C_\sigma} \simeq \rr^{|\sigma|}_{\geq 0}\,$
and extended by the gluing relations 
$ \,\overline{C}_\sigma \,\,\, \cap \,\, \overline{C}_\tau 
\,\, =\,\, \overline{C}_{\sigma \,\cap\, \tau}$.
\end{proof}

We now turn to the tropical Grassmannian 
and prove our first main result.
We shall identify the simplicial complex
${\bf T}_n$ with the fan in 
Theorem \ref{bhsfan}.

\medskip

\noindent {\sl Proof of Theorem \ref{phylogenetic}: }
The Pl\"ucker ideal $I_{2,n}$ is generated by
the $\binom{n}{4} $ quadrics
$$  p_{ij}  p_{kl}\, \, - \,\, p_{ik}  p_{jl} \,\, + \,\,
p_{il}  p_{jk} \qquad
\hbox{for} \,\,\,
1 \leq i < j < k < l \leq n . $$
The tropicalization of this polynomial is the
disjunction of linear systems
\begin{eqnarray*}
&  w_{ij} +  w_{kl} \, \, = \,\, w_{ik} +  w_{jl} \,\, \leq \,\,
 w_{il}  +  w_{jk} \\
{\rm or} \,\, 
&  w_{ij} +  w_{kl} \, \, = \,\,  w_{il}  +  w_{jk} \,\, \leq \,\,
 w_{ik} +  w_{jl} \\
{\rm or} \,\, 
&  w_{ik} +  w_{jl}  \, \, = \,\,  w_{il}  +  w_{jk} \,\, \leq \,\,
w_{ij} +  w_{kl} .
\end{eqnarray*}
Every point $w$ on the
tropical Grassmannian $\,{\mathcal G}_{2,n} \,$
satisfies this for  all quadruples $i,j,k,l$, that is,
it satisfies   the four point condition  (\ref{fourpoint}).
The Additive Linkage Algorithm reconstructs
the unique semi-labeled tree $\sigma$ with $w \in C_\sigma$.
This proves that every relatively open cone 
of  $\,{\mathcal G}_{2,n}\,$ lies in 
the relative interior of a unique cone $C_\sigma$ of the fan 
${\bf T}_n$ in Theorem 
\ref{bhsfan}.

We need to prove that the fans ${\bf T}_n$ and
$\,{\mathcal G}_{2,n} \,$ are equal.
Equivalently, every cone $C_\sigma$ is actually a cone
in the Gr\"obner fan. This will be accomplished by
analyzing the corresponding initial ideal.
In view of (\ref{itsafan}), it suffices to
consider maximal faces $\sigma$ of ${\bf T}_n$.
Fix a trivalent tree $\sigma$ and a
weight vector $w \in C_\sigma$. Then, for every
quadruple $i,j,k,l$, the inequality in
(\ref{fourpoint}) is strict. This means
combinatorially that 
$\,\bigl\{\{i,l\},\{j,k\}\bigr\}\,$
is a  four-leaf subtree of $\sigma$.

Let $J_\sigma$ denote the ideal
by the quadratic binomials
$\, p_{ij} p_{kl} - p_{ik} p_{jl} \,$
corresponding to all  four-leaf 
subtrees of $\sigma$. Our discussion
shows that $\, J_\sigma \subseteq {\rm in}_w (I_{2,n}) $.
The proof will be complete by showing that
the two ideals agree:
\begin{eqnarray}
\label{treeinitial}
\, J_\sigma \,\, = \,\, {\rm in}_w (I_{2,n}) .
\end{eqnarray}
This identity will be proved by showing that
the two ideals have a common initial
monomial ideal, generated by
square-free quadratic monomials.

We may assume, without loss of generality, that
$-w$ is a strictly positive vector, corresponding
to a planar realization of the tree $\sigma$
in which the leaves $1,2,\ldots,n$ are
arranged in circular order to form a convex $n$-gon
(Figure 1).

\begin{figure}
\centerline{\scalebox{0.5}{\includegraphics{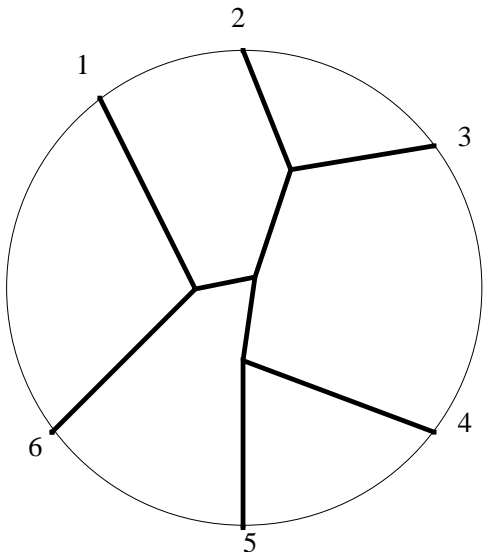}}}
\caption{A Circular Labeling of a Tree with Six Leaves}
\end{figure}

Let $M$ be the ideal generated by
the monomials $\,p_{ik} p_{jl} \, $ 
for $1 \leq i < j < k < l \leq n$.
These are the crossing pairs of edges in the $n$-gon.
By a classical construction of invariant theory,
known as \emph{Kempe's circular straightening law}
(see \cite[Theorem 3.7.3]{AIT}), there exists
a term order $\prec_{\rm circ}$ on $\zz[p]$ such that
\begin{equation}
\label{Mninitial}
M \quad = \quad {\rm in}_{\prec_{\rm circ}}(I_{2,n}). 
\end{equation}
Now, by our circular choice $w$ of realization of the tree $\sigma$,
the crossing monomials $\,p_{ik} p_{jl} \, $  appear
as terms in the binomial generators of $\, J_\sigma$.
Moreover, the term order $\prec_{\rm circ}$ on 
$\zz[p]$ refines the weight vector $w$. This implies
\begin{equation}
\label{chainofinc}
  {\rm in}_{\prec_{\rm circ}} ( {\rm in}_w (I_{2,n}))
\quad = \quad
{\rm in}_{\prec_{\rm circ}} (I_{2,n})
\quad = \quad M \quad
 \subseteq \quad {\rm in}_{\prec_{\rm circ}} (J_\sigma).
\end{equation}
Using $\,J_\sigma \subseteq {\rm in}_w (I_{2,n}) \,$
we conclude that equality holds in 
(\ref{chainofinc}) and in (\ref{treeinitial}).
\qed

\bigskip

The simplicial complex $\Delta(M)$ represented by
the squarefree monomial ideal $M$ is 
an iterated cone over the boundary of the 
polar dual of the \emph{associahedron};
see \cite[page 132]{AIT}. The facets of $\Delta(M)$
are the triangulations of the $n$-gon.
Their number is the common degree of the
ideals $I_{2,n}$, $J_\sigma$ and $M$:
$$ \hbox{the $(n-2)^{\rm nd}$ Catalan number}
 \quad = \quad \frac{1}{n-1} \binom{2n-4}{n-2}. $$
The reduced Gr\"obner basis
of (\ref{Mninitial}) has come to recent prominence
as a key example in the
Fomin-Zelevinsky theory of \emph{cluster algebras}
\cite{FZ}. Note also:

\smallskip

\begin{cor} There exists a maximal cone
in the Gr\"obner fan of
the Pl\"ucker ideal $I_{2,n}$ which contains,
up to symmetry, all cones of
$\,{\mathcal G}_{2,n}$.
\end{cor}

\begin{proof}
The cone corresponding to
the initial ideal  (\ref{Mninitial})
has this property.
\end{proof}

\begin{cor}
\label{inibinoisprime}
Every initial binomial ideal
of $I_{2,n}$ is a prime ideal.
\end{cor}

\begin{proof}
If ${\rm in}_w(I_{2,n})$ is a  binomial ideal
then $w$ must satisfy the four point conditions
(\ref{fourpoint}) with strict inequalities.
Hence $\,{\rm in}_w(I_{2,n}) \, = \, J_\sigma \,$
for some semi-labeled trivalent tree $\sigma$.
The ideal $J_\sigma$ is radical and equidimensional
because its initial ideal
$\ M \, = \, {\rm in}_{\prec_{\rm circ}}(J_\sigma)\,$
is radical and equidimensional (unmixed).

To show that $J_\sigma$ is  prime, we proceed as follows.
For each edge $e$ of the tree $\sigma$
we introduce an indeterminate $y_e$. Consider
the polynomial ring
$$\zz[y] \quad = \quad \zz \bigl[\,y_e \,:\, e \,\,\hbox{edge of}\,\,
\sigma \,\bigr]. $$
Let $\psi$ denote the homomorphism
$\zz[p] \rightarrow \zz[y]$ which sends
$p_{ij}$ to the product of all indeterminates
$y_e$ corresponding to edges on the
unique path between leaf $i$ and leaf $j$.
We claim  that $\, {\rm kernel}(\psi) \,= \,J_\sigma $.

A direct combinatorial argument shows that
the  convex polytope corresponding to the
toric ideal $\, {\rm kernel}(\psi) \,$
has a canonical triangulation 
into $\,\frac{1}{n-1} \binom{2n-4}{n-2}\,$
unit simplices (namely, $\Delta(M)$).
Hence $ \,{\rm kernel}(\psi) \,$
and  $\,J_\sigma \,$ are both unmixed of the
same dimension and the same degree.
Since $\, {\rm kernel}(\psi) \,$
is obviously contained in $\,J_\sigma $,
it follows that the two ideals are equal.
\end{proof}

\begin{cor}
\label{charfree}
The tropical Grassmannian 
$\,{\mathcal G}_{2,n}\,$
is characteristic-free.
\end{cor}

This means that we can consider
the Pl\"ucker ideal $I_{2,n}$ in the
polynomial ring $K[p]$ over any
ground field $K$ when computing its
tropical variety. All generators
$\, p_{ij} p_{kl} - p_{ik} p_{jl} \,$
of the initial binomial ideals
$\, J_\sigma \,$ have coefficients
$+1$ and $-1$, so $J_\sigma \otimes k$ contains
no monomial in $k[p]$, even if
$\,{\rm char}(k) > 0$.

\section{The Grassmannian of 3-planes in 6-space}

In this section we study the case $d=3$ and $n=6$.
The Pl\"ucker ideal $I_{3,6}$ is minimally 
generated by $35$ quadrics in the polynomial ring
in $20$ variables,
$$ \zz[p] \quad = \quad \zz[p_{123}, p_{124}, \ldots,p_{456}]. $$
We are interested in the $10$-dimensional fan 
${\mathcal G}_{3,6}$ which consists of all vectors
$w \in \rr^{20}$ such that
${\rm in}_w(I_{3,6})$ is monomial-free. The
four-dimensional quotient fan ${\mathcal G}_{3,6}''$ 
sits in $\,\rr^{20}/{\rm image}(\phi) \simeq \rr^{14}$ and 
is a fan over the three-dimensional polyhedral 
complex  ${\mathcal G}_{3,6}'''$. Our aim is to
prove Theorem \ref{G36}, which states that $\,{\mathcal G}_{3,6}''' \,$
consists of $65$ vertices,
$550$ edges, $1395$ triangles and $1035$ tetrahedra.

We begin by listing the  vertices.
Let $E$ denote the set of $20$ standard basis
vectors $e_{ijk}$ in $\rr^{\binom{6}{3}}$. For each
$4$-subset $\{i,j,k,l\}$ of $\{1,2,\ldots,6 \}$ we set
$$ f_{ijkl} \quad = \quad 
e_{ijk} \, + \,e_{ijl} \, + \, e_{ikl} \, + \, e_{jkl} . $$
Let $F$ denote the set of these $15$ vectors.
Finally consider any of the $15$ \emph{tripartitions}
$\,\{ \{i_1,i_2\}, \{i_3,i_4\}, \{i_5,i_6 \}\} \,$
of $\,\{1,2,\ldots,6\}$ 
and define the vectors
\begin{eqnarray*} &  g_{i_1 i_2 i_3 i_4 i_5 i_6}
\,\, := \,\, f_{i_1 i_2 i_3 i_4} 
\,+ \, e_{i_3 i_4 i_5}  \,+ \, e_{i_3 i_4 i_6}  \\
\ \hbox{and} \quad & 
 g_{i_1 i_2 i_5 i_6 i_3 i_4}
\,\, := \,\, f_{i_1 i_2 i_5 i_6} 
\,+ \, e_{i_3 i_5 i_6}  \,+ \, e_{i_4 i_5 i_6}.
\end{eqnarray*}
This gives us another set $G$ of $30$ vectors.
All $65$ vectors in $\,E \,\cup \, F \,\cup \, G \,$
are regarded as elements of the quotient space
$\, \rr^{\binom{6}{3}}/{\rm image}(\phi) \,\simeq \,\rr^{14}$.
Note that
$$ g_{i_1 i_2 i_3 i_4 i_5 i_6} \quad = \quad
 g_{ i_3 i_4 i_5 i_6 i_1 i_2} \quad = \quad g_{i_5 i_6 i_1 i_2 i_3 i_4 }. $$
Later on, the following identity will turn out to be important in the proof of Theorem~\ref{G36Detail}:
\begin{equation}
\label{notewid}
 g_{i_1 i_2 i_3 i_4 i_5 i_6} \, + \,  g_{i_1 i_2 i_5 i_6 i_3 i_4}
\quad = \quad
 f_{i_1 i_2 i_3 i_4}  \, + \,
f_{i_1 i_2 i_5 i_6}  \, + \,
f_{i_3 i_4 i_5 i_6} . 
\end{equation}
Lemma \ref{sixtyfive} and other results
in this section were found by computation.

\begin{lemma}
\label{sixtyfive}
The set of vertices of ${\mathcal G}_{3,6}$ equals
$\,E \,\cup \, F \,\cup \, G $.
\end{lemma}

We next describe all the $550$ edges of 
the tropical Grassmannian $\, {\mathcal G}_{3,6} $.
\begin{itemize}
\item[(EE)]
There are $90$ edges like $\{e_{123}, e_{145}\}$
and $10$ edges like $\{e_{123},e_{456}\}$,
for a total of $100$ edges connecting pairs of vertices
both of which are in $E$. (By the word ``like'', we will always mean ``in the $S_6$ orbit of, where $S_6$ permutes the indices $\{ 1, 2, \ldots 6 \}$.)
\item[(FF)] This class consists of $45$ edges like
$\,\{f_{1234}, f_{1256}\}$.
\item[(GG)] Each of the $15$ tripartitions gives exactly one edge, like
$\, \{g_{123456}, g_{125634} \}$.
\item[(EF)] There are $60$ edges like $\,\{ e_{123}, f_{1234} \}\,$
and $60$ edges like $\,\{ e_{123}, f_{1456}\} $, for a total
of $120$ edges connecting a vertex in $E$ to a vertex in $F$.
\item[(EG)] This class consists of $180$ edges like
$\,\{e_{123}, g_{123456}\}$. The intersections of the index triple of the $e$ vertex with the three index pairs of the $g$ vertex must have cardinalities $(2,1,0)$ \emph{in this cyclic order}.
\item[(FG)] This class consists of $90$ edges like $\, \{f_{1234}, g_{123456} \}$.
\end{itemize}

\begin{lemma}
The $1$-skeleton of ${\mathcal G}_{3,6}'''$ is the graph
with the $550$ edges above.
\end{lemma}

Let $\Delta$ denote the \emph{flag complex} specified by
the graph in the previous lemma. Thus $\Delta $ is
the simplicial complex on $\, E \cup F \cup G\,$ whose
faces are subsets $\sigma$ with the property that each 
$2$-element subset of $\sigma$ is 
one of the $550$ edges. We will see that
${\mathcal G}_{3,6}$ is a subcomplex 
homotopy equivalent to $\Delta$.

\begin{lemma} The flag complex $\Delta$ has
$1,410$ triangles, $1,065$ tetrahedra,
$15$ four-dimensional simplices, and it has no faces
of dimension five or more.
\end{lemma}

The facets of $\Delta$ are grouped into seven symmetry classes:

\smallskip

\noindent {\sl Facet FFFGG:}
There are $15$ four-dimensional simplices, one for each 
partition of $\{1,\ldots,6\}$ into three pairs.
An example of such a  tripartition is  $\{\{1,2\},\{3,4\},\{5,6\}\}$. 
It gives the facet
 $\, \{ f_{1234},f_{1256} ,f_{3456} ,g_{123456}, g_{125634} \}$.
 The $75$ tetrahedra contained in these $15$
four-simplices are not facets of $\Delta$.

\smallskip

The remaining $990$ tetrahedra in $\Delta$ are facets
and they come in six classes:

\smallskip

\noindent {\sl Facet EEEE:}
There are $30$ tetrahedra like
$\,\{e_{123}, e_{145}, e_{246}, e_{356} \}$.

\noindent {\sl Facet EEFF1:}
There are $90$ tetrahedra like
$\, \{ e_{123}, e_{456}, f_{1234}, f_{3456} \}$.

\noindent {\sl Facet EEFF2:}
There are $90$ tetrahedra like
$\,\{ e_{125}, e_{345}, f_{3456}, f_{1256} \}$.

\noindent {\sl Facet EFFG:}
There are $180$ tetrahedra like
$ \, \{ e_{345}, f_{1256}, f_{3456}, g_{123456} \}$.

\noindent {\sl Facet EEEG:} There are
$240$ tetrahedra like
$\,\{ e_{126}, e_{134}, e_{356},  g_{125634} \}$.

\noindent {\sl Facet EEFG: }
There are $360$ tetrahedra like
$\, \{ e_{234}, e_{125}, f_{1256}, g_{125634} \}$.

\smallskip

While $\Delta$ is an abstract simplicial complex on the vertices of $\mathcal{G}_{3,6}'''$, it is not embedded as a simplicial complex because relation 
(\ref{notewid}) shows that the five vertices of the four dimensional simplices only span three dimensional space. Specifically, they form a bipyramid with the 
F-vertices as the base and the G-vertices as the two cone points.

We now modify the flag complex $\Delta$ to a new 
simplicial complex $\Delta'$  which has pure dimension three and reflects the situation described in the last paragraph.
The complex $\Delta'$ is obtained from $\Delta$
by removing the $15$ FFF-triangles  $\, \{ f_{1234},f_{1256} ,f_{3456} \}$,
along with the $30$ tetrahedra FFFG and 
the $15$ four-dimensional facets FFFGG containing the FFF-triangles. In $\Delta'$, the bipyramids are each divided into three tetrahedra arranged around the GG-edges.
The following theorem implies both Theorem \ref{G36}
and Corollary \ref{G36notflag}.

\begin{thm} \label{G36Detail}
The tropical Grassmannian ${\mathcal G}_{3,6}'''$ equals
the simplicial complex $\Delta'$. It is not a flag complex
because of the $15$ missing FFF-triangles. The homology of ${\mathcal G}_{3,6}'''$ is concentrated in (top) dimension 3; $H_3({\mathcal G}_{3,6}''', \zz)=\zz^{126}$.
\end{thm}

The integral homology groups were computed independently by Michael Joswig 
and Volkmer Welker. We are grateful for their help.

This theorem is proved by an explicit computation. 
The correctness of the result can be verified by the following method.
One first
checks that the seven types of cones described above
are indeed Gr\"obner cones of $I_{3,6}$
whose initial ideals are monomial-free. Next
one checks that the list is complete. 
This relies on a result in
\cite{EKLW} which guarantees that ${\mathcal G}_{3,6}$ 
is connected in codimension $1$. The completeness check is done
by computing the link of each of the known 
classes of triangles. Algebraically, this amounts
to computing the (truly zero-dimensional) tropical
variety of ${\rm in}_w(I_{3,6})$ where $w$ is any point
in the relative interior of the triangular cone in question.
For all but one class of triangles the link consists of three points,
and each neighboring $3$-cell is found to be already among
our seven classes. The links of the triangles are as follows:

\smallskip

\noindent {\sl Triangle EEE: } The link  of
$\{e_{146}, e_{256}, e_{345}\}$  consists of 
$\,e_{123}$, $ g_{163425}$, $g_{142635}$.

\noindent {\sl Triangle EEF:} The link  of
$\{e_{256}, e_{346}, f_{1346}\}$ consists of
$\,f_{1256}$, $ g_{132546}$, $ g_{142536} $.

\noindent {\sl Triangle EEG: } The link  of
$\{e_{156}, e_{236}, g_{142356} \}$ consists of
$e_{124}$, $ e_{134}$, $ f_{1456}$.

\noindent {\sl Triangle EFF: } The link  of
$\{e_{135}, f_{1345}, f_{2346}\}$ consists of 
$\,e_{236}$, $ e_{246}$, $g_{153426}$.

\noindent {\sl Triangle EFG: } The link  of
$\{e_{235}, f_{2356}, g_{143526}\}$ consists of 
$e_{145}$, $f_{1246}$,  $e_{134}$.

\noindent {\sl Triangle FFG: } The link  of
$ \{f_{1236}, f_{1345}, g_{134526} \}$ consists of
$e_{126}$, $e_{236}$,  $g_{132645}$.

\noindent {\sl Triangle FGG: } The link  of
$ \{f_{1456}, g_{142356}, g_{145623}\}$ consists of
$f_{2356}$ and  $f_{1234}$.

\smallskip

The FGG triangle lies in the interior 
of our bipyramid FFFGG
and is incident to two of the three
FFGG tetrahedra which make up the
triangulation of that bipyramid. 
It is not contained in any other facet
of ${\mathcal G}_{3,6}'''$.

The $15$ bipyramids are
responsible for various counterexamples
regarding ${\mathcal G}_{3,6}$. This
includes the failure of
Corollaries \ref{isflagc} and
\ref{inibinoisprime} to hold
for $d \geq 3$.

\begin{prop}
\label{radicalnotprime}
Not every initial binomial ideal of $I_{3,6}$ is prime.
More precisely, if $w$ is any vector in the relative
interior of an FFGG cone then $in_w(I_{3,6})$ is
the intersection of two distinct codimension $10$ 
primes in  $\zz[p]$.
\end{prop}

\begin{proof}
We may assume that $\,w =  f_{1256} +  f_{3456} +  g_{123456} + 
g_{125634}$. Explicit computation (using 
\cite[Corollary 1.9]{GB+CP}) reveals
that $\,{\rm in}_w(I_{3,6}) \,$ is generated by
\begin{eqnarray*}
&
  p_{124} p_{135} - p_{123} p_{145},   \,\,\,
p_{123} p_{146} - p_{124} p_{136},   \,\,\,
   p_{125} p_{136} - p_{126} p_{135}, \\
&  p_{125} p_{146} - p_{126} p_{145},   \,\,\,
 p_{135} p_{146} - p_{136} p_{145},   \,\,\,
p_{123} p_{245} - p_{124} p_{235}, \\ &
p_{123} p_{246} - p_{124} p_{236},   \,\,\,
p_{126} p_{235} - p_{125} p_{236},   \,\,\,
p_{125} p_{246} - p_{126} p_{245},   \\ &
p_{134} p_{235} - p_{135} p_{234},   \,\,\,
p_{136} p_{234} - p_{134} p_{236},   \,\,\,
p_{136} p_{235} - p_{135} p_{236},  \\ &
p_{134} p_{245} - p_{145} p_{234},   \,\,\,
p_{134} p_{246} - p_{146} p_{234},   \,\,\,
p_{146} p_{245} - p_{145} p_{246},   \\ &
p_{135} p_{346} - p_{136} p_{345},   \,\,\,
p_{146} p_{345} - p_{145} p_{346},   \,\,\,
p_{135} p_{245} - p_{145} p_{235},  \\ &
p_{135} p_{256} - p_{156} p_{235},   \,\,\,
p_{156} p_{245} - p_{145} p_{256},   \,\,\,
p_{135} p_{456} - p_{145} p_{356},   \\ &
p_{136} p_{246} - p_{146} p_{236},   \,\,\,
p_{136} p_{256} - p_{156} p_{236},   \,\,\,
p_{146} p_{256} - p_{156} p_{246},  \\ &
p_{136} p_{456} - p_{146} p_{356},   \,\,\,
p_{235} p_{246} - p_{236} p_{245},   \,\,\,
p_{235} p_{346} - p_{236} p_{345},   \\ &
p_{245} p_{346} - p_{246} p_{345},   \,\,\,
p_{235} p_{456} - p_{245} p_{356},   \,\,\,
p_{246} p_{356} - p_{236} p_{456}, \\ &
p_{136} p_{245} - p_{135} p_{246},   \,\,\,
p_{145} p_{236} - p_{135} p_{246},   \,\,\,
p_{146} p_{235} - p_{135} p_{246},   \\ &
p_{123} p_{456} - p_{124} p_{356}\quad \hbox{and} \quad
p_{134} p_{256} - p_{156} p_{234}.
\end{eqnarray*}
The ideal $\,{\rm in}_w(I_{3,6})\,$ is the intersection of the two 
codimension $10$ primes
\begin{eqnarray*}
P  \quad = \quad
 {\rm in}_w (I_{3,6}) \,\,\,\, + & \!\!\!\!
\langle \,p_{125} p_{346} - p_{126} p_{345} \,\rangle 
\phantom{dadadada}  {\rm and}
\\
Q \quad = \quad
 {\rm in}_w (I_{3,6}) \,\,\,\, + &
\langle \,p_{135}, p_{136}, p_{145}, p_{146}, 
 p_{235}, p_{236}, p_{245}, p_{246}\, \rangle .
\end{eqnarray*}
The degrees of the ideals
$P$, $ Q$ and $I_{3,6}$ are
$\,38$, $4$ and $42$ respectively.
\end{proof}

We close this section with one more counterexample
arising from the triangulated bipyramid in ${\mathcal G}_{3,6}'''$.
It was proved in  \cite{BZ} 
that the $d \times d$-minors of a generic $d \times n$-matrix
form a \emph{universal Gr\"obner basis}.
A question left open in that paper is whether the
maximal minors also form a \emph{universal sagbi basis}.
It is well-known that they form
a sagbi basis for a specific term order.
See  \cite[Theorem 11.8]{GB+CP} and the
discussion  in Section 6 below. The question was whether
the sagbi basis property holds for all other
term orders. We show that the answer is ``no'':
the maximal minors are not a universal sagbi basis.

\begin{cor}
\label{notuniversalsagbi}
There exists a term order on $18$ variables such that 
the $3 \times 3$-minors of a generic
$3 \times 6$-matrix are not a 
sagbi basis in this term order.
\end{cor}

\begin{proof}
Consider the $3 \times 6$-matrix 
in \cite[Example 1.8 and Proposition 3.13]{SZ}:
$$ W \quad = \quad \begin{pmatrix}
2 & 1 & 2 & 1 & 0 & 0 \\
1 & 2 & 0 & 0 & 2 & 1 \\
0 & 0 & 1 & 2 & 1 & 2 
\end{pmatrix} $$
Let $ w \in \rr^{\binom{6}{3}}$ be its
vector of \emph{tropical $3\times 3$-minors}.
The coordinates of $ w$ are
\begin{eqnarray*}
 w_{ijk} \quad = \,\,\,
& {\rm min}
\bigl\{\,
  W_{1i}+W_{2j}+W_{3k},\,
 W_{1i}+W_{3j}+W_{2k},\,
 W_{2i}+W_{1j}+W_{3k}, \\ & \qquad
 W_{2i}+W_{3j}+W_{1k},\,
 W_{3i}+W_{1j}+W_{2k},\,
 W_{3i}+W_{2j}+W_{1k}
\bigr\}. 
\end{eqnarray*}
This vector represents the centroid of our bipyramid: 
$\,\, w \,= \, g_{123456} \,\, + \,\, g_{125634}$.
We consider the $3 \times 3$-minors of the following
matrix of indeterminates:
\begin{equation}
\label{XYZmatrix}
\begin{pmatrix}
\, x_1 & x_2 & x_3 & x_4 & x_5 &  x_6 \, \\
\,y_1 & y_2 & y_3 & y_4 & y_5 &  y_6 \,\\
\,z_1 & z_2 & z_3 & z_4 & z_5 &  z_6 
\end{pmatrix}
\end{equation}
The initial forms of its $3 \times 3$-minors
with respect to the weights $W$ are
\begin{eqnarray*} &
p_{123} =  z_1 x_2 y_3 \, , \,\,\, 
p_{124} =  z_1 x_2 y_4 \, , \,\,\, 
p_{125} =  y_1 z_2 x_5 \, , \,\,\, 
p_{126} =  y_1 z_2 x_6 \, , \\ &
p_{134} = -z_1 y_3 x_4 \, , \,\,\, 
p_{135} = -z_1 y_3 x_5 \, , \,\,\, 
p_{136} = -z_1 y_3 x_6 \, , \,\,\, 
p_{145} = -z_1 y_4 x_5 \, , \\ &
p_{146} = -z_1 y_4 x_6 \, , \,\,\, 
p_{156} =  z_1 x_5 y_6 \, , \,\,\, 
p_{234} = -z_2 y_3 x_4  \, , \,\,\, 
p_{235} = -z_2 y_3 x_5 \, , \\ &
p_{236} = -z_2 y_3 x_6 \, , \,\,\, 
p_{245} = -z_2 y_4 x_5 \, , \,\,\, 
p_{246} = -z_2 y_4 x_6 \, , \,\,\, 
p_{256} =  z_2 x_5 y_6 \, , \\ &
p_{345} = -z_3 y_4 x_5 \, , \,\,\, 
p_{346} = -z_3 y_4 x_6 \, , \,\,\, 
p_{356} =  y_3 z_5 x_6 \, , \,\,\, 
p_{456} =  y_4 z_5 x_6.
\end{eqnarray*}
Fix any term order $\prec$ which refines $W$.
The criterion in \cite[\S 11]{GB+CP} will
show that the $3 \times 3$-minors are not
a  sagbi basis  with respect to $\prec$.
The toric ideal of algebraic relations on 
the twenty monomials above
is precisely the prime $P$ in the proof
of Proposition \ref{radicalnotprime}.
The ideal $P$ strictly contains
${\rm in}_w(I_{3,6})$. Both have codimension $10$
but their degrees differ by $4$.
Using \cite[Theorem 11.4]{GB+CP} we conclude
that the $3 \times 3$-minors are not a sagbi basis for $\prec$.
\end{proof}

\section{Tropical Planes}

The Grassmannian $G_{d,n}$ is the parameter space
for all $d$-dimensional linear planes in $K^n$.
We now prove the analogous statement in
tropical geometry (Theorem \ref{grassmain}). 
But there are also crucial differences between 
the classical planes and tropical planes.
For instance, most  tropical planes are not 
complete intersections of tropical hyperplanes
(see Example \ref{funnytree} and Proposition \ref{caterpillar}).
Our combinatorial theory of tropical $d$-planes
is a direct generalization of the
Buneman representation of trees
(the $d=2$ case) and thus offers
mathematical tools for possible future
applications in phylogenetic analysis.

\medskip

\noindent {\sl Proof of Theorem \ref{grassmain}. }
The tropical Grassmannian  $\, {\mathcal G}_{d,n}'\,$
is a fan of dimension $(n-d)d$
in $\,\rr^{\binom{n}{d}} /\rr (1,1,\ldots,1) \,\simeq\,
\rr^{\binom{n}{d}-1}$.
We begin by describing the map
which takes a point $w $  in  $\, {\mathcal G}_{d,n}'\,$
to the associated tropical $d$-plane $L_w \subset \rr^n$.
Given $w$, we consider the tropical polynomials
\begin{equation}
\label{tropcirc}
 F_J(x_1,\ldots,x_n) \quad = \quad 
\sum_{j \in J}  w_{J \setminus \{ j \} } \cdot x_j ,
\end{equation}
where $J$ runs over all subsets of
cardinality $d+1$ in $[n]$. We define $L_w$ 
as the subset of $\rr^n$ which is
the intersection of the $\binom{n}{d+1}$ 
tropical hypersurfaces $\,{\mathcal T}(F_J)$.
We claim that $L_w$ is a tropical $d$-plane.
Pick a point $\xi \in (K^*)^{\binom{n}{d}}$ which is a zero
of $I_{d,n}$ and satisfies $w  = {\rm deg}(\xi)$. The
$d$-plane $X$ represented by $\xi$ is cut out by the
$\binom{n}{d+1}$  linear equations  derived from Cramer's rule:
\begin{equation}
\label{honestcirc}
f_J(x_1,\ldots,x_n) \quad = \quad
\sum_{j \in J}  \,\pm \,\xi_{J \setminus \{ j \}} \cdot x_j  \quad = \quad 0
\end{equation}
The tropicalization of this linear form is the
tropical polynomial in (\ref{tropcirc}), in symbols,
${\rm trop}(f_J) =  F_J$.  It is known that
the $f_J$ form a universal Gr\"obner basis
for the ideal they generate
\cite[Proposition 1.6]{GB+CP}.  Therefore,
Corollary \ref{FiTG} 
shows that
$L_w$ is indeed a tropical
$d$-plane. In fact, we have
$$ {\rm deg}(X) \quad = \quad L_w \quad = \quad
L_{{\rm deg}(\xi)}. $$
This proves that the map $w \mapsto L_w$
surjects the tropical Grassmannian
onto the set of all tropical $d$-planes,
and it is the only such map which is
compatible with the classical bijection between
$G_{d,n}$ and the set of $d$-planes in $K^n$.

It remains to be shown that 
the map $w \mapsto L_w$ is injective.
We do this by constructing the inverse map.
Suppose we are given  $L_w$ as a subset of $\rr^n$. We need
to reconstruct the coordinates $w_{i_1 \cdots i_d}$
of $w$ up to a global additive constant.
Equivalently, for any $(d-1)$-subset
$I$ of $[n]$ and any pair 
$j, k \in [n] \backslash I$, we need to 
reconstruct the real number
$\, w_{I \cup \{j\}} - w_{I \cup \{k\}}$.

Fix a very large positive rational number $M $
and consider the $(n-d+1)$-dimensional plane
defined by $\, x_i = M \,$ for $i \in I$.
The intersection of this plane
with $L_w$ contains at least one point $x \in \rr^n$,
and this point
can be chosen to satisfy
$\,x_j \ll M \,$ for all $j \in [n] \backslash I$.
This can be seen by solving the
$d-1$ equations  $\,x_i = t^M \,$
on any  $d$-plane $X \subset K^n$
which tropicalizes to $L_w$.

Now consider
the tropical polynomial (\ref{tropcirc})
with $\,J\, =\, I \,\cup \, \{j,k\}$.
Since $x$ lies ${\mathcal T}(F_J)$, and since
$\,{\rm max}(x_j,x_k) \ll M = x_i \,$ for all $i \in I$,
we conclude
$$  w_{J \backslash \{k\}} \, + \,  x_k \quad = \quad
w_{J \backslash \{j\}} \, + \,  x_j . $$
This shows that the desired differences can be read off from 
the point $x$:
\begin{equation}
\label{differences}
\, w_{I \cup \{j\}} - w_{I \cup \{k\}} \quad = \quad
x_j \,\, - \,\, x_k. 
\end{equation}
We thus reconstruct
$\, w \in {\mathcal G}_{d,n}\,$ by locating
$\binom{n}{d-1}$ special points on  $\,L_w$.
\qed

\medskip 

The above proof offers an (inefficient) algorithm for
computing the map $w \mapsto L_w$, namely, by
intersecting all $\binom{n}{d+1}$ tropical
hypersurfaces ${\mathcal T}(F_J)$.
Consider the case $d=2$. Here
the $\binom{n}{3}$ tropical polynomials 
$F_J$ in  (\ref{tropcirc}) are
$$ F_{ijk} \quad = \quad
 w_{ij} \cdot x_k \,\, + \,\,
 w_{ik} \cdot x_j \,\, + \,\,
 w_{jk} \cdot x_i .$$
The tropical hypersurface 
$\,{\mathcal T}(F_{ijk}) \,$
is the solution set to the
linear system
\begin{eqnarray*}
 & w_{ij}\, + \, x_k \,\,=\,\,
w_{ik}\, + \, x_j \,\,\leq\,\,
w_{jk}\, + \, x_i  \\
\quad \hbox{or} \quad
 & w_{ij}\, + \, x_k \,\,= \,\,
w_{jk}\, + \, x_i \,\, \leq\,\,
w_{ik}\, + \, x_j  \\
\quad \hbox{or} \quad
 & 
w_{ij}\, + \, x_k \,\,=\,\,
w_{jk}\, + \, x_i  \,\,\leq\,\,
w_{ik}\, + \, x_j.
\end{eqnarray*}
The conjunction of these $\binom{n}{3}$ linear
systems can be solved efficiently by a variant
of the tree reconstruction algorithm in
\cite{Bu}. If $r$ and $s \in \RR^n / \RR (1, \ldots, 1)$ are 
vertices of this tree connected by an edge $e$ then 
$r=s+c \sum_{i \in S} e_i$ for some $c>o$ where 
$S \subset [n]$ is the set of leaves on the $s$ side of $e$.
We regard the tree as a metric space by 
assigning the length $c$ to edge $e$.
The length of each edge is measured in
 lattice distance, so we get the tree with
 metric $-2w$.  

\begin{cor} \label{reallygetatree}
Let $w$ be a point in ${\mathcal G}_{2,n}$
which lies in the cone $C_\sigma$ for
some tree $\sigma$. The image 
of $L_w$ in $\rr^n/\rr (1,\ldots,1)$ is a
tree  of combinatorial type $\sigma$. 
\end{cor}

The bijection $w \mapsto L_w $ 
of Theorem \ref{grassmain} is a  higher-dimensional 
generalization of recovering a phylogenetic tree from
pairwise distances among $n$ leaves.  For instance,
for $d=3$ we can think of $w$ as data giving
a proximity measure for any triple among $n$ ``leaves''.
The image of $L_w$ in $\rr^n/\rr (1,\ldots,1)$ is a
``phylogenetic surface'' which is a geometric
representation of such data.

The tropical Grassmannians ${\mathcal G}_{d,n}$
and ${\mathcal G}_{n-d,n}$ are isomorphic 
because the ideals $I_{d,n}$  and $I_{n-d,n}$
are the same after signed complementation of
Pl\"ucker coordinates. Theorem 
\ref{grassmain} allows us to define the
\emph{dual} $(n-d)$-plane  $L^*$ of a tropical
$d$-plane $L$ in $\rr^n$. 
If $L = L_w$ then $L^* = L_{w^*}$ where
$w^*$ is the vector whose
$([n]\backslash I)$-coordinate is the $I$-coordinate
of $w$, for all $d$-subsets $I$ of $[n]$. One can check that a tropical hyperplane $\sum a_i \cdot x_i=0$ contains $L^*$ iff $(a_i) \in L_w$ and that $(L^*)^*=L$.

\begin{ex}
\label{funnytree}
Let $w = e_{12} + e_{34} + e_{56}$ in $\rr^{\binom{6}{2}}$.
Then $L_w$ is a tropical $2$-plane in $\rr^6$. Its image
in $\rr^6/\rr (1,\ldots,1)$ is a tree as in Figure 1, of type
$\,\sigma \, = \,\bigl\{ \{12,3456\},   \{34,1256\},   \{56,1234\}\bigr\}$.
The  Pl\"ucker vector dual to $w$ is
$$ w^* \quad = \quad e_{3456} \, + \, 
e_{1256} \, + \, e_{1234} \,\quad \in \,\, \, {\mathcal G}_{4,6}
\,\,\subset\, \, \rr^{\binom{6}{4}}. $$
We shall compute the tropical $4$-plane 
$L_{w^*}$ by applying the algorithm in the
proof of Theorem \ref{grassmain}. There
are $6$ tropical polynomials $F_J$ as in 
(\ref{tropcirc}), namely,
\begin{eqnarray*}
F_{12345} \quad = &   0 \cdot x_1 \, + \,
 0 \cdot x_2 \, + \,  0 \cdot x_3 \, + \, 0 \cdot x_4 \, + \,  1 \cdot x_5 \\
F_{12346} \quad = &   0 \cdot x_1 \, + \,
 0 \cdot x_2 \, + \,  0 \cdot x_3 \, + \, 0 \cdot x_4 \, + \,  1 \cdot x_6 \\
F_{12356} \quad = &   0 \cdot x_1 \, + \,
 0 \cdot x_2 \, + \,  1 \cdot x_3 \, + \, 0 \cdot x_5 \, + \,  0 \cdot x_6 \\
F_{12456} \quad = &   0 \cdot x_1 \, + \,
 0 \cdot x_2 \, + \,  1 \cdot x_4 \, + \, 0 \cdot x_5 \, + \,  0 \cdot x_6 \\
F_{13456} \quad = &   1 \cdot x_1 \, + \,
 0 \cdot x_3 \, + \,  0 \cdot x_4 \, + \, 0 \cdot x_5 \, + \,  0 \cdot x_6 \\
F_{23456} \quad = &   1 \cdot x_2 \, + \,
 0 \cdot x_3 \, + \,  0 \cdot x_4 \, + \, 0 \cdot x_5 \, + \,  0 \cdot x_6 
\end{eqnarray*}
The tropical $4$-plane $L_{w^*}$  
is the intersection of these six tropical hyperplanes:
$$ 
{\mathcal T}(F_{12345})\,\cap \,
{\mathcal T}(F_{12346})\,\cap \,
{\mathcal T}(F_{12356})\,\cap \,
{\mathcal T}(F_{12456})\,\cap \,
{\mathcal T}(F_{13456})\,\cap \,
{\mathcal T}(F_{23456}).
$$
We claim that $\, L_{w^*}\,$ \emph{is not
a complete intersection}, i.e.,
there do no exist two tropical linear
forms $F$ and $F'$ such that
$\, L_{w^*} \, = \, {\mathcal T}(F) \, \cap \,{\mathcal T}(F')$.
A tropical linear form  $\, F  = a_1 x_1 + \cdots + a_6 x_6\,$
vanishes on the dual $4$-plane $L_{w^*}$ if and only if
the point
$a = (a_1,\ldots,a_6)$ lies in the $2$-plane $L_w$.
There are $9$ types of such tropical linear forms
$F$, one for each of the $9$ edges of the tree $L_w$.
For instance, the bounded edge
$\{56,1234\}$ represents the tropical forms
$$\, F \,\, = \,\,
\alpha \cdot (x_1 + x_2 + x_3 + x_4) \,
 + \, \beta \cdot (x_5 + x_6) 
 \quad \hbox{where} \,\,\,
0 < \alpha \leq \beta. $$
By checking all pairs of the $9$ edges, we find
that any conceivable intersection
$\,{\mathcal T}(F) \, \cap \,{\mathcal T}(F')\,$
must contain a $5$-dimensional cone like
$\{x_1 +c = x_2   \ll x_3,x_4,x_5,x_6\}$,
$\{x_3 +c = x_4  \ll x_1,x_2,x_5,x_6\}$ or
$\{x_5 +c = x_6  \ll x_1,x_2,x_3,x_4\}$.
\end{ex}

This example can be generalized as follows.

\begin{prop} 
\label{caterpillar}
Let $L_w$ be a tropical $\, 2$-plane
in $\rr^n$ whose tree is not
combinatorially isomorphic to  $\,\sigma = 
\bigl\{ \{1,\ldots,i\},\{i+1,\ldots,n\} \, :\,
i = 2,3,\ldots,n-2 \bigr\}$. Then
the dual tropical $(n-2)$-plane $L_{w^*}$ is not
a complete intersection.
\end{prop}

The special tree $\sigma$ in Proposition \ref{caterpillar}
is called the \emph{caterpillar} in the phylogenetic literature
(see Figure 2). 

\begin{proof}
Suppose for contradiction that $L_{w^*}$ is the intersection of the 
hyperplanes $\sum a_i \cdot x_i =0$ and $\sum b_i \cdot x_i=0$. 
The vectors $a = (a_1,\ldots,a_n)$ and
$b= (b_1,\ldots,b_n)$, regarded as
elements of $ \RR^n/\RR(1,\ldots,1)$, lie in the
tree $L_w$. Denote by $\gamma$ the path through $L_w$ from $a$ to $b$.
Since $L_w$ is not a caterpillar tree,
the path $\gamma$ goes through fewer than $n-1$ edges, so deleting those edges divides $L$ into fewer then $n$ connected components. Thus, there are two leafs of $L$, call them $j$ and $k$, such that the none of the edges of $\gamma$ 
separate $j$ from $k$. 
Every edge of $\gamma$ connects two points $r$ and $s$ with $s=r+c \sum_{i \in S} e_i$ where, in each case, either $j$ and $k$ both lie in $S$ or neither do. Thus, $a_j-a_k=b_j-b_k$. Therefore, the intersection of the hyperplanes $\sum a_i \cdot x_i=0$ and $\sum b_i \cdot x_i=0$ contains every point $(x_i)$ with $x_j+a_j=x_k+a_k$ and $x_i-x_j$ sufficiently positive for all $i \neq j, k$. But this is a codimension one subset of $\RR^n/\RR(1,\ldots,1)$ and we know that $L_{w^*}$ is pure of codimension two.
\end{proof}

\begin{figure}
\centerline{\scalebox{0.7}{\includegraphics{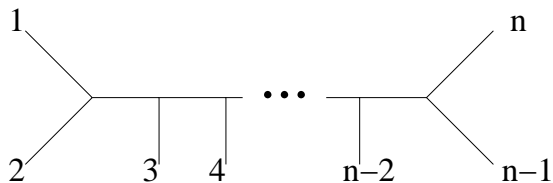}}}
\caption{A Caterpillar Tree}
\end{figure}

Our next goal is to give a combinatorial encoding
of tropical planes. The basic object in our combinatorial encoding
is a $d$-partition
$\,\{A_1,\ldots,A_d\}$. By a
\emph{$d$-partition} we mean
an  unordered partition
of $[n]$ into $d$ subsets $A_i$.
Let $L_w$ be a tropical $d$-plane
and $F$ a maximal cell of $L_w$.
Thus $F$ is a $d$-dimensional
convex polyhedron in $\rr^n$.
The affine span of $F$ is a
$d$-dimensional  affine space which is defined by
equations of the special form
$$  x_k \,\, - \,\,  x_j \quad = \quad
w_{J \backslash \{j\}} \,\, - \,\,
  w_{J \backslash \{k\}} \qquad \hbox{(the right
hand side is a constant)} $$
Such a  system of equations defines a $d$-partition
$\{A_1,\ldots,A_d\}$, namely, 
two indices $j$ and $k$ lie in the
same block $A_i$ if and only
if the difference $x_k - x_j$ is constant on $F$.
The number of blocks clearly equals
$d$, the dimension of $F$.

\begin{rmk}
\label{boundedfaces}
A maximal face $F$ of $L_w$ is uniquely specified
by its $d$-partition $\{A_1,\ldots,A_d\}$.
It is a (bounded) polytope in $\rr^n$ if and only
if $\,|A_i| \geq 2 \,$ for all $i$.
Hence a tropical $d$-plane $L_w \subset \rr^n$
has no bounded $d$-faces if $n \leq 2d-1$.
\end{rmk}

We define the   \emph{type} of a tropical $d$-plane $L$,
denoted ${\rm type}(L)$ to be the set of
all $d$-partitions  arising from the maximal faces of $L$.
If $d=2$ and $L=L_w$ with $w \in  C_\sigma$ then
$ {\rm type}(L) $ is precisely the set 
$\sigma $ together with the
pairs $ \{ \{i\}, [n] \backslash \{i\} \}$
representing the unbounded edges of the tree $L$.
This follows from Corollary \ref{reallygetatree}.
Thus ${\rm type}(L)$ generalizes the
Buneman representation of semi-labeled trees
(Section 4) to higher-dimensional tropical planes $L$.

The type of a tropical plane $L_w$ is a strong combinatorial invariant, 
but it does not uniquely determine the cone of ${\cal G}_{d,n}$
which has $w$ in its relative interior. 
We will see this phenomenon in the example below.

\vskip .1cm

\begin{ex}
\label{threeofseven}
We present three of the seven types in
${\mathcal G}_{3,6}$. In each case we display
${\rm type}(L_w)$ with the
$15$ obvious tripartitions $\,\bigl\{i,j,[6] \backslash \{i,j\}\bigr\}$
removed.

We begin with a type which we call the \emph{sagbi type}:
\begin{eqnarray*}
 &  \bigl\{ 
 \{1 , 23 , 456 \},  \{1 , 56 , 234 \}, 
 \{2 , 13 , 456 \},  \{2 , 56 , 134 \}, \\
\hbox{EEFF1}  : & \,\,
 \{3 , 12 , 456 \},  \{3 , 56 , 124 \}, 
 \{4 , 12 , 356 \},  \{4 , 56 , 123 \}, \\ & \!\!
 \{5 , 12 , 346 \},  \{5 , 46 , 123 \}, 
 \{6 , 12 , 345 \},  \{6 , 45 , 123\}, \{12 , 34 , 56 \}  \bigr\}
\end{eqnarray*}
The next type is the \emph{bipyramid type}.
All three tetrahedra in a bipyramid {\sl FFFGG}
have the same type listed below. As the faces of $\mathcal{G}_{3,6}'''$ contain those $w$ inducing different initial ideals ${\rm in}_w (I_{d,n})$, 
this example demonstrates that ${\rm type}(L_w)$ does not determine 
${\rm in}_w (I_{d,n})$.

\begin{eqnarray*}
 &  \bigl\{ 
 \{1 , 34 , 256\},  \{1 , 56 , 234\}, 
 \{2 , 34 , 156\},  \{2 , 56 , 134\} \\
\hbox{FFGG}:  & \,\,\,\,
 \{3 , 12 , 456\},  \{3 , 56 , 124\},
 \{4 , 12 , 356\},  \{4 , 56 , 123\} \\ &
 \{5 , 12 , 346\},  \{5 , 34 , 126\},
 \{6 , 12 , 345\},  \{6 , 34 , 125\},
 \{12 , 34 , 56 \}  \bigr\}
\end{eqnarray*}
For all but one of the seven types in ${\mathcal G}_{3,6}$,
the tropical plane $L_w$ has $28$ facets.
The only exception is the type {\sl EEEE}.
Here the tropical plane $L_w$ has only $27$
facets, all of them unbounded.
\begin{eqnarray*}
 &  \bigl\{ 
\{ 1, 23,456 \},  \{ 1,234,56 \},
\{ 2,13,456 \},  \{ 2,135,46 \} \\
\hbox{EEEE} : \qquad \quad & \,\,\,\,
\{ 3,12,456 \},  \{ 3,126,45 \}, 
\{ 4,26,135 \},  \{ 4,126,35 \} \\ & \quad \,
\{ 5,16,234 \},  \{ 5,126,34 \}, 
\{ 6,15,234 \},  \{6, 135, 24\} \bigr\}
\end{eqnarray*}
\end{ex}

\section{Dependence on the characteristic}

Our definition of the tropical Grassmannian implicitly depended on the fields
$K$ and $k$. The ideal $I_{d,n}$ makes sense over any field and has the same generators (the classical Pl\"ucker relations). Nonetheless, the properties of the initial ideal $in_w(I_{d,n})$ might depend on $k$, in particular, whether or not this ideal contains a monomial might depend on the charcteristic of $k$. 
Hence, whether or not $w \in {\mathcal G}_{d,n}$ might depend on the characteristic of $k$.

In Corollary \ref{charfree} we saw that 
this does not happen for $d = 2$,
and it follows from the explicit
computations in Section 6 that this
does not happen for ${\mathcal G}_{3,6}$
either. In both of these cases,  the tropical Grassmannian
is  characteristic-free. Another result that we observed
in both of these nice cases is that it was enough to look
at quadratic polynomials in $I_{d,n}$ to define
the tropical Grassmannian. We shall see below that
the same results do not hold for the next case
${\mathcal G}_{3,7}$. We summarize our result this in
the following theorem.

\begin{thm}
Let $ d= 2$ or $d = 3$ and $n = 3$. Then
every monomial-free initial ideal
of $I_{d,n}$ is generated by quadrics,
and the tropical Grassmannian
${\mathcal G}_{d,n}$ is characteristic-free.
Both of these properties fail for
$d \geq 3$ and $n \geq 7 $.
\end{thm}

\begin{proof}
It suffices to consider the case $ d  = 3$ and $n = 7$.
An easy lifting argument will extend our example
to the general case $d \geq 3$ and $n \geq 7 $.
The Pl\"ucker ideal $I_{3,7}$ is minimally
generated by $140$ quadrics in a polynomial
ring $k[ p_{123},p_{124}, \ldots, p_{567}]$
in $35$ unknowns over an arbitrary field $k$.

We fix the following zero-one vector. The appearing
triples are gotten gotten by a cyclic shift, and they correspond
to the lines in the \emph{Fano plane}:
$$ w \quad = \quad e_{124} + e_{235} + e_{346} + e_{457} + e_{156} + e_{267} + e_{137}
\quad \in \,\,\, \rr^{\binom{6}{3}} . $$
We next compute the initial ideal $\, {\rm in}_w(I_{3,7})$
under the assumption that the characteristic
of $k$  is zero. In a computer algebra system,
this is done by  computing the reduced Gr\"obner basis of 
$I_{3,7}$ over the field of rational numbers
with respect to the (reverse lexicographically
refined) weight order defined by $-w$.
The reduced Gr\"obner basis is found to have precisely $196$ elements,
namely,  $140 $ quadrics, $52$ cubics, and $4$ quartics.
The initial ideal $\, {\rm in}_w(I_{3,7})\,$ is generated by
the $w$-leading forms of the $196$ elements in that
Gr\"obner basis.

Among the $52$ cubics in the Gr\"obner basis of $I_{3,7}$,
we find the special cubic
\begin{eqnarray*} & f \quad =  \quad
 2 \cdot p_{123} p_{467} p_{567}
\,-\, p_{367} p_{567} \underline{p_{124}}
- p_{167} p_{467} \underline{p_{235}}
- p_{127} p_{567} \underline{p_{346}}  \\ & \quad
- p_{126} p_{367} \underline{p_{457}}
- p_{237} p_{467} \underline{p_{156}}
+ p_{134} p_{567} \underline{p_{267}}
+ p_{246} p_{567} \underline{p_{137}}
 \,\, + \, p_{136} \underline{p_{267}} \underline{p_{457}}.
\end{eqnarray*}
Since $\,{\rm char}(k) \not= 2$, the leading
form of this polynomial is the monomial
$$ {\rm in}_w(f) \quad = \quad p_{123} p_{467} p_{567}. $$
This proves that $w$ is not in the tropical 
Grassmannian  ${\mathcal G}_{3,7}$.

On the other hand, suppose now that
the characteristic of $k$ equals two. In that case,
the leading form of $f$ is a polynomial with seven terms
$$ {\rm in}_w(f) \quad = \quad
\,-\, p_{367} p_{567} \underline{p_{124}}
- p_{167} p_{467} \underline{p_{235}} -
\, \cdots \,
+ p_{246} p_{567} \underline{p_{137}}.
$$
This is not a monomial. In fact, none of the leading forms
of the $196$ Gr\"obner basis elements is a monomial.
This proves that the initial ideal  ${\rm in}_w(I_{3,7})$
contains no monomial, or equivalently,
that  $w$  lies in the tropical Grassmannian
 ${\mathcal G}_{3,7}$ when ${\rm char}(k) = 2$.
In fact, there is no inclusion in either direction between the
tropical Grassmannians ${\mathcal G}_{3,7}$
in characteristic two and in characteristic zero.
To see this, we modify our vector $w$ as follows:
$$ w' \quad = \quad w - e_{124} 
 \quad = \quad
 e_{235} + e_{346} + e_{457} + e_{156} + e_{267} + e_{137}
\quad \in \,\,\, \rr^{\binom{6}{3}} . $$
Then  $\, in_{w'}(f) =  2 \cdot p_{123} p_{467} p_{567}  - p_{367} p_{567} \underline{p_{124}}$,
which is not monomial if ${\rm char}(k) = 0$, but it is a 
monomial if ${\rm char}(k) = 2$. This shows that $w'$ does not  lie in
${\mathcal G}(3,7)$ if the characteristic of $k$ is two.
By recomputing the Gr\"obner basis in characteristic zero, we find that the initial 
ideal $\,in_{w'}(I_{3,7})\,$ contains no monomial, and hence 
does  lie in ${\mathcal G}(3,7)$ if the characteristic of $k$ is zero.

The above argument also shows that, in any
characteristic, either ${\rm in}_w(I_{3,6})$ or
${\rm in}_{w'}(I_{3,6})$ will be a monomial-free
initial ideal which has a minimal generator
of degree three. Quadrics do not suffice
for $d \geq 3$ and $n \geq 7$.
\end{proof}

It is worth taking a moment to think about the intuitive geometry behind this argument. Let $B$ be any subset of $\binom{[n]}{d}$; we can study the collection of points on the Grassmannian $G(d,n)$ over the field $k$ where the Pl\"ucker coordinate $P_I$ is nonzero exactly for those $I \in B$. Such points exist exactly if $B$ is the set of bases of a matroid of rank $d$ on $n$ points realizable over $k$. 

Thus, when the characteristic of $k$ is $2$ there is a point 
$x \in G_{3,7}$ with $x_{ijk}=0$ exactly when $i$, $j$ and $k$ are 
collinear in the Fano plane and no such point should 
exist in characterisitc other than $2$. Passing to the tropicalization,
one would expect that in characteristic $2$ there should be a point $y \in {\mathcal G}_{d,n}$ with $y_{ijk}=\infty$ for $i$, $j$ and $k$ collinear in the Fano plane and $y_{ijk}=0$ otherwise. Intuitively, $w$ is a perturbation of $y$ so that $w_{ijk}$ is 1 and not $\infty$.

\end{document}